\documentclass[12pt]{article}
\usepackage{amssymb}

\usepackage{amsmath,amssymb,amsfonts}
\usepackage{amsmath,amssymb,amsfonts}

\newtheorem{theorem}{Theorem}
\newtheorem{proposition}[theorem]{Proposition}

\begin{document}

\author{Marcos Salvai\thanks{%
Partially supported by \textsc{fonc}y\textsc{t}, Antorchas, \textsc{%
ciem\thinspace (conicet)} and \textsc{sec}y\textsc{t\thinspace (unc)}.}}
\title{On the geometry of the space of oriented lines of the hyperbolic space}
\date{}
\maketitle

\begin{abstract}
Let $H$ be the $n$-dimensional hyperbolic space of constant sectional
curvature $-1$ and let $G$ be the identity component of the isometry group
of $H$. We find all the $G$-invariant pseudo-Riemannian metrics on the space
$\mathcal{G}_{n}$ of oriented geodesics of $H$ (modulo orientation
preserving reparametrizations). We characterize the null, time- and
space-like curves, providing a relationship between the geometries of $%
\mathcal{G}_{n}$ and $H$. Moreover, we show that $\mathcal{G}_{3}$ is
K\"{a}hler and find an orthogonal almost complex structure on $\mathcal{G}%
_{7}$.
\end{abstract}

\noindent%

\smallskip%

\noindent%

\bigskip

\begin{center}
\textbf{1. The space of geodesics of a Hadamard manifold}
\end{center}

\bigskip

\noindent Let $M$ be a Hadamard manifold (a complete simply connected
Riemannian manifold with nonpositive sectional curvature) of dimension $n+1$%
. An \emph{oriented geodesic} $c$ of $M$ is a complete connected totally
geodesic oriented submanifold of $M$ of dimension one. We may think of $c$
as the equivalence class of unit speed geodesics $\gamma :\mathbf{R}%
\rightarrow M$ with image $c$ such that $\left\{ \dot{\gamma}\left( t\right)
\right\} $ is a positive basis of $T_{\gamma \left( t\right) }c$ for all $t$%
. Let $\mathcal{G}=\mathcal{G}\left( M\right) $ denote the space of all
oriented geodesics of $M$. The space of geodesics of a manifold all of whose
geodesics are periodic with the same length is studied with detail in \cite
{besse}. The geometry of the space of oriented lines of Euclidean space is
studied in \cite{gk, salvaimm, salvai2}.

Let $T^{1}M$ be the unit tangent bundle of $M$ and $\xi $ the spray of $M$,
that is, the vector field on $T^{1}M$ defined by $\xi \left( v\right)
=\left. d/dt\right| _{0}\gamma _{v}^{\prime }\left( t\right) $, where $%
\gamma _{v}$ is the unique geodesic in $M$ with initial velocity $v$.
Clearly, $\mathcal{G}$ may be identified with the set of oriented leaves of
the foliation of $T^{1}M$ induced by $\xi $. By \cite{keilhauer}, if $M$ is
Hadamard, this foliation is regular in the sense of Palais \cite{palais}.
Hence, $\mathcal{G}$ admits a unique differentiable structure of dimension $%
2n$ such that the natural projection $T^{1}M\rightarrow \mathcal{G}$ is a
submersion.

Fix $o\in M$ and let Exp $:T_{o}M\rightarrow M$ denote the geodesic
exponential map. Let $S=\left\{ v\in T_{o}M\mid \left\| v\right\| =1\right\}
\cong S^{n}$. We identify as usual $T_{v}S\cong v^{\bot }\subset T_{o}M$.
Hence, $TS\cong \left\{ \left( v,x\right) \mid v\in S\text{ and }%
\left\langle v,x\right\rangle =0\right\} $. Let $F:TS\rightarrow \mathcal{G}$
be defined by
\[
F\left( v,x\right) =\left[ \gamma \right] \text{,}
\]
where $\gamma $ is the unique geodesic in $M$ with initial velocity $\tau
_{0}^{1}v$ (here $\tau $ denotes parallel transport along the geodesic $%
t\rightarrow $ Exp\thinspace $\left( tx\right) $ of $M$). This is called the
minitwistor construction in \cite{hitchin}. Keilhauer proved in \cite
{keilhauer} that $F$ is a diffeomorphism.

\bigskip

\begin{center}
\textbf{2. The geometry of $\mathcal{G}$ for the hyperbolic space}
\end{center}

\bigskip

\noindent Let $H=H^{n+1}$ be the hyperbolic space of constant sectional
curvature $-1$ and dimension $n+1$. Consider on $\Bbb{R}^{n+2}$ the basis $%
\left\{ e_{0},e_{1},\dots ,e_{n+1}\right\} $ and the inner product whose
associated norm is given by $\left\| x\right\| =\left\langle
x,x\right\rangle =-x_{0}^{2}+x_{1}^{2}+\dots +x_{n+1}^{2}$. Then $H=\left\{
x\in \Bbb{R}^{n+2}\mid \left\| x\right\| =-1\text{ and }x_{0}>0\right\} $
with the induced metric. Let $G$ be the identity component of the isometry
group of $H$, that is,
\[
G=O_{o}\left( 1,n+1\right) =\left\{ g\in O\left( 1,n+1\right) \mid \left(
ge_{0}\right) _{0}>0\text{ and }\det g>0\right\} \text{.}
\]
In the following we denote $\mathcal{G}_{m}=\mathcal{G}\left( H^{m}\right) $
(or simply $\mathcal{G}$ if no confusion is possible). The group $G$ acts on
$\mathcal{G}$ as follows: $g\left[ \gamma \right] =\left[ g\circ \gamma
\right] $. This action is transitive, since $H$ is two-point homogeneous,
and smooth, since $G$ acts smoothly on $T^{1}H$.

Let $\gamma _{o}$ be the geodesic in $H$ with $\gamma _{o}\left( 0\right)
=e_{0}$ and initial velocity $e_{1}\in T_{e_{0}}H$. The isotropy subgroup of
$G$ at $c_{o}:=\left[ \gamma _{o}\right] $ is
\[
G_{o}=\left\{ \text{diag}\left( T_{t},A\right) \mid t\in \Bbb{R},A\in
SO_{n}\right\} \cong \Bbb{R}\times SO_{n}\text{,}
\]
where $T_{t}=\left(
\begin{array}{ll}
\cosh t & \sinh t \\
\sinh t & \cosh t
\end{array}
\right) $. Therefore we may identify $\mathcal{G}$ with $G/G_{o}$ in the
usual way. Let $\frak{g}$ be the Lie algebra of $G$ and let
\[
\frak{g}_{o}=\left\{ \text{diag}\left( tR,A\right) \mid t\in \Bbb{R},A\in
so_{n}\right\}
\]
be the Lie algebra of $G_{o}$ (here $R=\left(
\begin{array}{ll}
0 & 1 \\
1 & 0
\end{array}
\right) $). Let $B$ be the bilinear form on $\frak{g}$ defined by $B\left(
X,Y\right) =\frac{1}{2}\,$tr\thinspace $\left( XY\right) $, which is
well-known to be a multiple of the Killing form of $\frak{g}$, hence
nondegenerate. Besides, the canonical projection $\pi :G\rightarrow H$, $\pi
\left( g\right) =g\left( e_{0}\right) $, is a pseudo-Riemannian submersion.

Let $\frak{g}=\frak{g}_{0}\oplus \frak{h}$ be the orthogonal decomposition
with respect to $B$. Then
\[
T_{c_{o}}\mathcal{G}=\frak{h}:=\left\{ x_{h}+y_{v}\mid x,y\in \Bbb{R}%
^{n}\right\} ,
\]
where for column vectors $x,y\in \Bbb{R}^{n}$,
\[
x_{h}=\left(
\begin{array}{cc}
0_{2} & \left( x,0\right) ^{t} \\
\left( x,0\right) & 0_{n}
\end{array}
\right) \text{ \ \ and\ \ \ \ }y_{v}=\left(
\begin{array}{cc}
0_{2} & \left( 0,y\right) ^{t} \\
\left( 0,-y\right) & 0_{n}
\end{array}
\right)
\]
(here the exponent $t$ denotes transpose and $0_{m}$ the $m\times m$ zero
matrix). We chose this notation since $x_{h}$ and $y_{v}$ are horizontal and
vertical, respectively, tangent vectors in $T_{\left( e_{0},e_{1}\right)
}\left( T^{1}H\right) $ with respect to the canonical projection $%
T^{1}H\rightarrow H$.

\begin{theorem}
\label{homogeneous}For each $n\geq 1$ there exists a $G$-invariant
pseudo-Riemannian metric $g_{1}$ on $\mathcal{G}_{n+1}$ whose associated
norm at $c_{o}$ is given by
\[
\left\| x_{h}+y_{v}\right\| _{1}=\left| x\right| ^{2}-\left| y\right| ^{2}%
\text{.}
\]
For $n=2,$ if one identifies $\Bbb{R}^{2}=\Bbb{C}$ as usual, there exists a $%
G$-invariant metric $g_{0}$ on $\mathcal{G}_{3}$ whose associated norm at $%
c_{o}$ is given by
\[
\left\| x_{h}+y_{v}\right\| _{0}=\left\langle ix,y\right\rangle \text{.}
\]
For $n\ne 2$, any $G$-invariant pseudo-Riemannian metric on $\mathcal{G}%
_{n+1}$ is homothetic to $g_{1}$. Any $G$-invariant pseudo-Riemannian metric
on $\mathcal{G}_{3}$ is of the form $\lambda g_{0}+\mu g_{1}$ for some $%
\lambda ,\mu \in \Bbb{R}$ not simultaneously zero.

All the metrics are symmetric and have split signature $\left( n,n\right) $.
In particular, $\mathcal{G}$ does not admit any $G$-invariant Riemannian
metric and the geodesics in $\mathcal{G}$ through $c_{o}$ are exactly the
curves $s\mapsto \exp _{G}\left( sX\right) c_{o}$, for $X\in \frak{h}$.
\end{theorem}

\noindent%
\textbf{Proof. }One computes easily that $B\left( X,X\right) =\left\|
X\right\| _{1}$ for all $X\in \frak{h}$. Since $B$ is $G$-invariant, $g_{1}$
defines a $G$-invariant metric on $\mathcal{G}$.

Let $Z=$ diag\thinspace $\left( R,0_{n}\right) $, $\frak{m}=\{ \text{diag}%
\left( 0_{2},A\right)\mid A\in so_{n}\} $ and $\frak{g}_{\lambda }=\{ U\in
\frak{g}\mid \text{ad}_{Z}U=\lambda U\} $. One verifies that $\frak{g}_{0}=%
\frak{g}_{o}$ and $\frak{g}_{\pm 1}=\left\{ x_{h}\pm x_{v}\mid x\in \Bbb{R}%
^{n}\right\} $. Moreover, one has the decompositions
\[
\frak{g}_{0}=\Bbb{R}Z\oplus \frak{m}\;\;\;\text{and\ \ \ }\frak{h}=\frak{g}%
_{1}\oplus \frak{g}_{-1}\text{,}
\]
which are preserved by the the action of $\frak{m}$. Hence $\frak{h}$ is $%
\frak{g}_{0}$-invariant.

Since $B$ is nondegenerate and $G_{o}$ is connected, any other
pseudo-Rie\-mannian metric $g$ on $\mathcal{G}$ has the form $g\left(
U,V\right) =B\left( TU,V\right) $ for some $T:\emph{h}\rightarrow \emph{h}$
commuting with ad$_{Z}$ and ad$_{\frak{m}}$. In particular, $T$ preserves $%
\frak{g}_{\pm 1}$. We call $T_{\pm }$ the restrictions of $T$ to the
corresponding subspaces. Under the identification $\frak{g}_{\pm 1}\equiv
\Bbb{R}^{n}$, $x_{h}\pm x_{v}\equiv x$, the action of $\frak{m}\equiv so_{n}$
on $\Bbb{R}^{n}$ is the canonical one. If $T_{\pm }\in Gl\left( \frak{g}%
_{\pm 1}\right) \equiv Gl\left( n,\Bbb{R}\right) $ commutes with every $A\in
so_{n}$, then either $T_{\pm }$ is a nonzero multiple of the identity or $%
n=2 $ and $T_{\pm }=a_{\pm }I_{2}+b_{\pm }J$ where $J=\left(
\begin{array}{cc}
0 & -1 \\
1 & 0
\end{array}
\right) $, for some not simultaneously zero constants $a_{\pm }$ and $b_{\pm
}$. Next we consider the case $n=2$ and show that $a_{+}=a_{-}$ and $%
b_{-}=-b_{+}$. For $x\ne 0$ we denote $x^{\pm }=x_{h}\pm x_{v}$ and compute
\begin{eqnarray*}
B\left( T\left( x^{+}\right) ,x^{-}\right) &=&B\left( \left(
a_{+}x+b_{+}ix\right) ^{+},x^{-}\right) \\
&=&a_{+}B\left( x^{+},x^{-}\right) +b_{+}B\left( \left( ix\right)
^{+},x^{-}\right) \\
&=&2a_{+}\left| x\right| ^{2}+0\text{.}
\end{eqnarray*}
Since $T$ must be symmetric with respect to $B$, this expression coincides
with $B\left( x^{+},T\left( x^{-}\right) \right) $, which by similar
computations equals $2a_{-}\left| x\right| ^{2}$. Hence $a_{+}=a_{-}$. Using
again the symmetry of $T$ in the case
\[
B\left( T\left( x^{-}\right) ,\left( ix\right) ^{+}\right) =B\left(
x^{-},T\left( ix\right) ^{+}\right)
\]
one obtains that $b_{-}=-b_{+}$. Finally, since $2\left( x_{h}+y_{v}\right)
=\left( x+y\right) ^{+}+\left( x-y\right) ^{-}$, one computes that the
metric associated with $T$ is homothetic to $g_{1}$ if $b_{+}=0$ and to $%
g_{o}$ if $a_{+}=0$. The case $n\ne 2$ is simpler since it does not involve $%
b_{\pm }$.

Next we show that for any of the metrics above, $\mathcal{G}$ is a symmetric
space. Let $G^{\uparrow }=\left\{ g\in O\left( 1,n+1\right) \mid \left(
ge_{0}\right) _{0}>0\right\} $ be the isometry group of $H$ and let $C=$
diag\thinspace $\left( I_{2},-I_{n}\right) \in G^{\uparrow }$, which induces
an involutive diffeomorphism $\widetilde{C}$ of $\mathcal{G}$ by $\widetilde{%
C}\left[ \gamma \right] =\left[ C\circ \gamma \right] $ fixing exactly $%
c_{o} $. If $n=2$, $C\in G$, hence $\widetilde{C}$ is clearly an isometry
for any $G$-invariant metric on $\mathcal{G}_{3}$. The same happens for $%
n\ne 2$. Indeed, in this case, up to homotheties, we have seen that the
unique metric on $\mathcal{G}_{m}$ with $m\ne 3$ comes from a multiple of
the Killing form of $\frak{g}$, which is invariant by the action of $%
G^{\uparrow }$. The statement regarding geodesics follows from the theory of
symmetric spaces, since conjugation by $C$ is an involutive automorphism of $%
\frak{g}$ whose $\left( -1\right) $-eigenspace is $\frak{g}_{0}$ and
preserves the given metrics. \hfill $\square $

\bigskip%

\noindent%
\textbf{Remarks. }a) In contrast with the space of oriented lines of $\Bbb{R}%
^{n}$, which only for $n=3,7$ admits pseudo-Riemannian metrics invariant by
the induced transitive action of a connected closed subgroup of the identity
component of the isometry group (see \cite{salvaimm}), $\mathcal{G}_{n}$
admits $G$-invariant metrics for all $n$.

\smallskip%

b) The metric $g_{0}$ is the analogue of the metric defined in the Euclidean
case in \cite{shepherd, gk1}. We will see below that also in the hyperbolic
case it admits a K\"{a}hler structure.

\smallskip%

c) For any complete simply connected Riemannian manifold $M$ of negative
curvature, the space $\mathcal{G}\left( M\right) $ of its oriented geodesics
has a canonical pseudo-Riemannian metric, which is in general only
continuous, see \cite{Kanai}. If $M$ is the hyperbolic space, then $g_{1}$
is the canonical metric on $\mathcal{G}$.

\smallskip%

d) If $H$ has dimension two, then $\mathcal{G}$ is isometric to the
two-dimensional de Sitter sphere.

\bigskip%

We recall some well-known facts about the imaginary border of the hyperbolic
space and the action of $G$ on it. For a geodesic $\gamma $ in $H$, $\gamma
\left( \infty \right) $ is defined to be the unique $z\in S^{n}$ such that $%
\lim_{t\rightarrow \infty }\gamma \left( t\right) /\gamma \left( t\right)
_{0}=e_{0}+z\in \Bbb{R}^{n+2}$. One defines analogously $\gamma \left(
-\infty \right) $. Sometimes we will identify $\Bbb{R}^{n+1}$ with $%
e_{0}^{\bot }$ and $S^{n}$ with $\left\{ e_{0}\right\} \times S^{n}$.

The group $G$ acts on $S^{n}$ by directly (that is, orientation preserving)
conformal diffeomorphisms. More precisely, any $g\in G$ induces the directly
conformal transformation $\widetilde{g}$ of $S^{n}$, well-defined by $%
\widetilde{g}\left( \gamma \left( \infty \right) \right) =\left( g\circ
\gamma \right) \left( \infty \right) $, and any directly conformal
transformation of $S^{n}$ can be realized in this manner.

\begin{proposition}
\label{proptoni}If $S$ is a subgroup of $G$ acting transitively on $\mathcal{%
G}$, then $S=G$.
\end{proposition}

\noindent%
\textbf{Proof. }By the main result of \cite{toni}, it suffices to show that $%
S$ acts irreducibly on $\Bbb{R}^{n+2}$. Suppose that $S$ leaves the
nontrivial subspace $V$ invariant. If $V$ is degenerate, then $V$ contains a
null line, say $\Bbb{R}\left( e_{0}+z\right) $, with $z\in \Bbb{R}^{n+1}$, $%
\left| z\right| =1$. Hence $S$ takes the oriented line $\left[ \gamma
\right] $ with $\gamma \left( \infty \right) =z$ to another line with the
same point at $\infty $. If $V$ is nondegenerate, either $V$ or its
complement (also $S$-invariant) intersects $H$. Let us call $%
H_{1}\varsubsetneq H$ the intersection, which is a totally geodesic
submanifold of $H$. Then $S$ takes any oriented line contained in $H_{1}$ to
a line contained in $H_{1}$. If $H_{1}$ is a point $p$, then $S$ takes any
line through $p$ to a line through $p$. Therefore the action of $S$ on $%
\mathcal{G}$ is not transitive. \hfill $\square $

\medskip%

\noindent%
\textbf{Remark. }The hyperbolic case contrasts with the Euclidean one: We
found in \cite{salvaimm} a pseudo-Riemannian metric on the space of oriented
lines of $\Bbb{R}^{7}=$ Im\thinspace $\Bbb{O}$ which is invariant by the
transitive action of $G_{2}\ltimes \Bbb{R}^{7}$, where $G_{2}$ is the
automorphism group of the octonions $\Bbb{O}$.

\bigskip

\begin{center}
\textbf{3. Null, space- and time-like curves}
\end{center}

\bigskip

\noindent In order to give a geometric interpretation for a curve in $%
\mathcal{G}$ endowed with some of the $G$-invariant metrics to be null,
space- or time-like, we introduce the following concept, which makes sense
for any Hadamard manifold.

\medskip%

\noindent%
\textbf{Definition.} Let $H$ be a Hadamard manifold. Given a smooth curve $c$
in $\mathcal{G}$ defined on the interval $I$, a function $\varphi :\mathbf{R}%
\times I\rightarrow H$ is said to be a \emph{standard presentation of }$c$
if $s\mapsto \alpha _{t}\left( s\right) :=\varphi \left( s,t\right) $ is a
unit speed geodesic of $H$ satisfying $c\left( t\right) =\left[ \alpha
_{t}\right] $ and $\left\langle \dot{\beta}\left( t\right) ,\dot{\alpha}%
_{t}\left( 0\right) \right\rangle =0$ for all $t\in I$, where $\beta \left(
t\right) =\varphi \left( 0,t\right) $.

\begin{proposition}
\label{presentation}Given a smooth curve $c:I\rightarrow \mathcal{G}$ and $p$
a point in the image of some \emph{(}any\emph{)} geodesic in the equivalence
class $c\left( t_{o}\right) $, there exists a standard presentation $\varphi
$ of $c$ such that $\varphi \left( 0,t_{o}\right) =p$.
\end{proposition}

\noindent%
\textbf{Proof.} Consider the submersion $\Pi :T^{1}H\rightarrow \mathcal{G}$%
, $\Pi \left( v\right) =[\gamma _{v}]$. Let $v\left( t\right) $ be a lift of
$c\left( t\right) $ to $T^{1}H$ with $v\left( t_{o}\right) \in T_{p}^{1}H$,
and let $\psi :\mathbf{R}\times I\rightarrow H$ be defined by $\psi \left(
s,t\right) =\gamma _{v\left( t\right) }\left( s\right) $. We look for a
function $f:I\rightarrow \mathbf{R}$ such that
\[
\varphi \left( s,t\right) =\psi \left( s+f\left( t\right) ,t\right)
\]
satisfies the required properties. Clearly $\alpha _{t}\left( s\right)
=\varphi \left( s,t\right) $ has unit speed and
\[
c\left( t\right) =\Pi \left( v\left( t\right) \right) =\Pi \left( \gamma
_{v\left( t\right) }^{\prime }\left( f\left( t\right) \right) \right)
=\left[ \alpha _{t}\right] \text{.}
\]
One can verify easily that taking as $f$ the solution of the differential
equation
\[
f^{\prime }\left( t\right) =-\frac{\left\langle \psi _{t}\left( f\left(
t\right) ,t\right) ,\psi _{s}\left( f\left( t\right) ,t\right) \right\rangle
}{\left\| \psi _{s}\left( f\left( t\right) ,t\right) \right\| ^{2}}
\]
(subindexes denote partial derivatives) with $f\left( t_{o}\right) =0$, then
$\varphi ( 0,t_{o}) =p$ and $\left\langle \dot{\beta}\left( t\right) ,\dot{%
\alpha}_{t}\left( 0\right) \right\rangle =0$ for all $t\in I$, where $\beta $
is as in the definition of the standard presentation. \hfill $\square $

\bigskip%

The following Proposition characterizes the null, time- and space-like
curves of $\mathcal{G}$, providing a relationship between the geometries of $%
\mathcal{G}$ and $H$.

\begin{proposition}
\label{relation}For the metric $g_{1}$, a smooth curve $c$ in $\mathcal{G}%
_{n}$ is null \emph{(}respectively, space-, time-like\emph{)} if and only
if, for any standard presentation, the rate of variation of the directions,
that is, $\left\| \frac{D}{dt}\dot{\alpha}_{t}\left( 0\right) \right\| $,
coincides with \emph{(}respectively, is smaller, larger than\emph{)} the
rate of displacement $\left\| \dot{\beta}\left( t\right) \right\| $ for all $%
t$ \emph{(}here $\frac{D}{dt}$ denotes covariant derivative along $\beta $%
\emph{)}.

For the metric $g_{0}$ on $\mathcal{G}_{3}$, a smooth curve $c$ in $\mathcal{%
G}_{3}$ is null \emph{(}respectively, space-\nolinebreak , time-like\emph{)}
if and only if, for any standard presentation,
\[
\left\{ \dot{\beta}\left( t\right) ,\frac{D}{dt}\dot{\alpha}_{t}\left(
0\right) ,\dot{\alpha}_{t}\left( 0\right) \right\}
\]
is linearly dependent \emph{(}respectively, positively, negatively oriented%
\emph{)} for all $t$.
\end{proposition}

\noindent%
\textbf{Proof. }Let $\left[ \gamma \right] $ be an oriented geodesic of a
Hadamard manifold and let $\mathcal{J}_{\gamma }$ be the space of Jacobi
fields along $\gamma $ orthogonal to $\dot{\gamma}$. First we show that $%
L_{\gamma }:\mathcal{J}_{\gamma }\rightarrow T_{\left[ \gamma \right] }%
\mathcal{G}$ given by
\begin{eqnarray}
L_{\gamma }\left( J\right) =\left( d/dt\right) _{0}\left[ \gamma _{t}\right]
\text{,}  \label{el}
\end{eqnarray}
where $\gamma _{t}$ is a variation of $\gamma $ by unit speed geodesics
associated with the Jacobi field $J$, is a well-defined vector space
isomorphism. Indeed, let $\mathcal{P}:T^{1}M\rightarrow \mathcal{G}$ be the
canonical projection, which is a smooth submersion, by definition of the
differentiable structure on $\mathcal{G}$. We compute
\[
\left( d/dt\right) _{0}\left[ \gamma _{t}\right] =\left( d/dt\right) _{0}%
\mathcal{P}\left( \dot{\gamma}_{t}\left( 0\right) \right) =d\mathcal{P}_{%
\dot{\gamma}\left( 0\right) }\left( \left( d/dt\right) _{0}\dot{\gamma}%
_{t}\left( 0\right) \right) .
\]
Now, let p $:T^{1}H\rightarrow H$ be the canonical projection and $\mathcal{K%
}:T_{\dot{\gamma}\left( 0\right) }\left( T^{1}H\right) \rightarrow \dot{%
\gamma}\left( 0\right) ^{\bot }\subset T_{\gamma \left( 0\right) }H$ the
connection operator. It is well-known that $\left( d\text{p},\mathcal{K}%
\right) :T_{\dot{\gamma}\left( 0\right) }\left( T^{1}H\right) \rightarrow
T_{\gamma \left( 0\right) }H\oplus \dot{\gamma}\left( 0\right) ^{\bot }$ is
a bijection and
\[
\left( d/dt\right) _{0}\dot{\gamma}_{t}\left( 0\right) =\left( d\text{p},%
\mathcal{K}\right) ^{-1}\left( J\left( 0\right) ,J^{\prime }\left( 0\right)
\right)
\]
(see for instance \cite{besse}). Therefore, $L_{\gamma }$ is well-defined.

Next we show that for any $J\in \mathcal{J}_{\gamma }$ one has
\begin{eqnarray}
\left\| L_{\gamma }\left( J\right) \right\| _{1} &=&\left\| J\left( 0\right)
\right\| ^{2}-\left\| J^{\prime }\left( 0\right) \right\| ^{2}  \label{nj} \\
\left\| L_{\gamma }\left( J\right) \right\| _{0} &=&\left\langle \dot{\gamma}%
\left( 0\right) \times J\left( 0\right) ,J^{\prime }\left( 0\right)
\right\rangle \text{.}  \nonumber
\end{eqnarray}
We may suppose without loss of generality that $c=c_{o}$ and $\gamma =\gamma
_{o}.$ Let $c^{\prime }\left( 0\right) =x_{h}+y_{v}$ with $x,y\in \Bbb{R}%
^{n} $. Then the Jacobi field along $\gamma _{o}$ satisfying $L_{\gamma
_{o}}\left( J\right) =c^{\prime }\left( 0\right) $ is the one determined by
\[
J\left( 0\right) =d\pi _{I}\left( x_{h}\right) \text{ and }J^{\prime }\left(
0\right) =d\pi _{I}\left( y_{h}\right) \text{,}
\]
where $\pi :G\rightarrow H$ is as before the canonical projection. In fact,
clearly, $\gamma _{t}\left( s\right) =\exp \left( tx_{h}\right) \exp \left(
ty_{v}\right) \gamma _{o}\left( s\right) $ is a variation of $\gamma _{o}$
by unit speed geodesics. Let us see that the associated Jacobi field is $J.$
Indeed,
\[
J\left( 0\right) =\left. \frac{d}{dt}\right| _{0}\gamma _{t}\left( 0\right)
=\left. \frac{d}{dt}\right| _{0}\exp \left( tx_{h}\right) e_{0}=d\pi
_{I}\left( x_{h}\right) \text{,}
\]
since $\gamma _{o}\left( 0\right) =e_{0}$, which is fixed by $\exp \left(
ty_{v}\right) $. If $\frac{D}{dt}$ denotes covariant derivative along $%
t\mapsto \gamma _{t}\left( 0\right) $ and $Z$ is as in the beginning of the
proof of Theorem \ref{homogeneous}, then
\begin{eqnarray*}
J^{\prime }\left( 0\right) &=&\left. \frac{D}{dt}\right| _{0}\dot{\gamma}%
_{t}\left( 0\right) =\left. \frac{D}{dt}\right| _{0}d\left( \exp \left(
tx_{h}\right) \exp \left( ty_{v}\right) \right) _{\pi \left( I\right) }e_{1}
\\
&=&\left. \frac{D}{dt}\right| _{0}d\exp \left( tx_{h}\right) d\pi _{I}\text{%
Ad}\left( \exp ty_{v}\right) Z \\
&=&d\pi _{I}\left. \frac{d}{dt}\right| _{0}e^{t\text{ad\thinspace }%
y_{v}}Z=d\pi _{I}\left[ y_{v},Z\right] =d\pi _{I}\left( y_{h}\right) \text{,}
\end{eqnarray*}
since $d\exp \left( tx_{h}\right) $ realizes the parallel transport and $%
d\pi _{I}\left( Z\right) =e_{1}$. Therefore (\ref{nj}) is true by Theorem
\ref{homogeneous}. Finally, suppose that $\varphi $ is a standard
presentation of $c$ and let $\alpha _{t},\beta $ be as above. Let $J_{t}$
denote the Jacobi field along $\alpha _{t}$ associated with the variation $%
\varphi $. Clearly, $\dot{c}\left( t\right) =L_{\alpha _{t}}\left(
J_{t}\right) $, $J_{t}\left( 0\right) =\frac{d}{dt}\varphi \left( 0,t\right)
=\dot{\beta}\left( t\right) $ and
\[
J_{t}^{\prime }\left( 0\right) =\left. \frac{D}{ds}\right| _{0}\frac{d}{dt}%
\varphi \left( s,t\right) =\frac{D}{dt}\left. \frac{d}{ds}\right|
_{0}\varphi \left( s,t\right) =\frac{D}{dt}\dot{\alpha}_{t}\left( 0\right) .
\]
Consequently, the proposition follows from (\ref{nj}). \hfill $\square $

\bigskip%

\noindent%
\textbf{A geometric invariant of }$\mathcal{G}$

\bigskip%

\noindent%
We have mentioned in the introduction that $\mathcal{G}\left( H^{n}\right) $
is diffeomorphic to $\Bbb{T}^{n}$, the space of all oriented lines of $\Bbb{R%
}^{n}$. For $n=3$ and $n=7$, we found in \cite{salvaimm} pseudo-Riemannian
metrics on $\Bbb{T}^{n}$ invariant by the induced transitive action of a
connected closed subgroup of $SO_{n}\ltimes \Bbb{R}^{n}$ (only for those
dimensions such metrics exist).

\begin{proposition}
For $n=3,7$, no metric on $\mathcal{G}_{n}$ invariant by the identity
component of the isometry group of $H^{n}$ is isometric to $\Bbb{T}^{n}$
endowed with any of the metrics above.
\end{proposition}

\noindent \textbf{Proof. }We compute now a pseudo-Riemannian invariant of $%
\mathcal{G}_{n}$ involving its periodic geodesics. For any $c\in \mathcal{G}$%
, let $A$ denote the subset of $T_{c}\mathcal{G}$ consisting of the
velocities of periodic geodesics of $\mathcal{G}$ though $c$. We show next
that the frontier of $A$ in $T_{c}\mathcal{G}$ is the union of two subspaces
of half the dimension of $\mathcal{G}$ intersecting only at zero. By
homogeneity we may suppose that $c=c_{o}$. Since by the proposition below $%
A=\left\{ \lambda x_{h}+x_{v}\mid x\in \Bbb{R}^{n}\text{, }\left| \lambda
\right| <1\right\} $, the frontier of $A$ is $\frak{g}_{1}\cup \frak{g}_{-1}$%
. On the other hand, we have computed in \cite{salvai2} that the analogue
invariant for $\Bbb{T}^{n}$ ($n=3,7$) is a subspace of half the dimension of
$\Bbb{T}^{n}$. Hence the proposition follows. \hfill $\square $

\medskip

\noindent \textbf{Remarks. }a) Of course we could have considered more
standard invariants, like the curvature or the isometry group, but we chose
this one since the geodesics can be described so easily.

b) Clearly the difference in the invariants is related to the fact that the
two horospheres through a point associated with opposite directions coincide
in the Euclidean case but are different in the hyperbolic case.

\begin{proposition}
A geodesic in $\mathcal{G}$ with initial velocity $x_{h}+y_{v}$ is periodic
if and only if $x=\lambda y$ for some $\lambda \in \Bbb{R}$ with $\left|
\lambda \right| <1$.
\end{proposition}

\noindent%
\textbf{Proof. }We may suppose that $x_{h}+y_{v}\ne 0$. We compute
that Ad\thinspace $\left( e^{tZ}\right) \left(
x_{h}+y_{v}\right)\allowbreak =x_{v}^{t}+y_{v}^{t}$, where
\[
x^{t}=\left( \cosh t\right) x+\left( \sinh t\right) y\text{ \ \ and \ \ }%
y^{t}=\left( \sinh t\right) x+\left( \cosh t\right) y\text{.}
\]
Now, there exists $s$ such that $\left\langle x^{s},y^{s}\right\rangle =0$
(take $\tanh \left( 2s\right) =-\frac{2\left\langle x,y\right\rangle}{\left|
x\right| ^{2}+\left| y\right| ^{2}} $). Hence $\left[
x_{h}^{s},y_{v}^{s}\right] =0$ and consequently
\[
\pi \exp \left( t\left( x_{h}^{s}+y_{v}^{s}\right) \right) =\pi \exp \left(
tx_{h}^{s}\right) \exp \left( ty_{v}^{s}\right) =\pi \exp \left(
tx_{h}^{s}\right) \text{,}
\]
which is a geodesic in $H$, in particular it is periodic only if it is
constant, or equivalently, only for $x^{s}=0$.

Since $Z\in \frak{g}_{0}$ and the metric is $G$-invariant, the geodesics
with initial velocities $x_{h}^{t}+y_{v}^{t}$ are simultaneously periodical
or not periodical for all $t$. Now, one verifies that $x^{s}=0$ if and only
if $x=\lambda y$ for some $\lambda \in \Bbb{R}$ with $\left| \lambda \right|
<1$ and the proposition follows. \hfill $\square $

\bigskip

\begin{center}
\textbf{4. Additional geometric structures on $\mathcal{G}$}
\end{center}

\bigskip

\noindent An almost Hermitian structure on a pseudo-Riemannian manifold $%
\left( M,g\right) $ is a smooth tensor field $J$ of type $\left( 1,1\right) $
on $M$ such that $J_{p}$ is an orthogonal transformation of $\left(
T_{p}M,g_{p}\right) $ and satisfies $J_{p}^{2}=-\,$id for all $p\in M$. If $%
\nabla $ is the Levi Civita connection of $\left( M,g\right) ,$ then $\left(
M,g,J\right) $ is said to be K\"{a}hler if $\nabla J=0$.

\bigskip

\noindent%
\textbf{A K\"{a}hler structure on }$\mathcal{G}\left( H^{3}\right) $

\bigskip

\noindent Let $\mathcal{G}=\mathcal{G}_{3}$ and let $j_{o}$ be the
endomorphism of $\frak{h}\equiv T_{c_{o}}\mathcal{G}\equiv \Bbb{C}\times
\Bbb{C}$ given by $j_{o}\left( z,w\right) =\left( iz,iw\right) $. One checks
that $j_{o}$ commutes with the action of $G_{o}$, is orthogonal for $g_{0}$
and $g_{1}$ and $j_{o}^{2}=-\,$id. Therefore $j_{o}$ defines an orthogonal
almost complex structure on $\mathcal{G}_{3}$ for any $G$-invariant metric
on it.

\begin{proposition}
The space $\left( \mathcal{G}_{3},J\right) $ is K\"{a}hler for any
pseudo-Riemannian $G$-invariant metric on $\mathcal{G}_{3}$.
\end{proposition}

\noindent \textbf{Proof. }We show that for every geodesic $\gamma $ in $%
\mathcal{G}_{3}$ and any parallel vector field $Y$ along $\gamma $, the
vector field $JY$ along $\gamma $ is parallel. By homogeneity we may suppose
that $\gamma \left( 0\right) =c_{o}$. Suppose that $\gamma \left( t\right)
=\exp \left( tX\right) c_{o}$ for some $X\in \frak{h}$. By a well-known
property of symmetric spaces, $Y=d\exp \left( tX\right) _{c_{o}}Y_{c_{o}}$.
Since $J$ is $G$-invariant, $JY=d\exp \left( tX\right) _{c_{o}}JY_{c_{o}}$
and thus $JY$ is parallel along $\gamma $, as desired. \hfill $\square $

\bigskip%

\noindent \textbf{An orthogonal almost complex structure on }$\mathcal{G}%
_{7} $

\bigskip%

\noindent We present another model of $\mathcal{G}_{n+1}$ endowed with the
metric $g_{1}$ and use it to define an orthogonal almost complex structure
on $\mathcal{G}_{7}$.

\medskip%

In the following we use the notations given before Proposition \ref{proptoni}
of concepts related to the imaginary border of $H$. We recall that $g\in G$
is called a transvection of $H$ if it preserves a geodesic $\gamma $ of $H$
and $dg$ realizes the parallel transport along $\gamma $, that is, $g\left(
\gamma \left( t\right) \right) =\gamma \left( t+s\right) $ for all $t$ and
some $s$ and $dg_{\gamma \left( t\right) }$ realizes the parallel transport
between $\gamma \left( t\right) $ and $\gamma \left( t+s\right) $ along $%
\gamma $. For any unit $v\in T_{e_{0}}H=e_{0}^{\bot }=\Bbb{R}^{n+1}$ the
transvections through $e_{0}\in H$ preserving the geodesic with initial
velocity $v$ form a one parameter subgroup $\phi _{t}$ such that the
corresponding one parameter group $\widetilde{\phi _{t}}$ of conformal
transformations of $S^{n}$ (which we also call transvections, by abuse of
notation) is the flow of the vector field on $S^{n}$ defined at $q\in S^{n}$
as the orthogonal projection of the constant vector field $v$ on $\Bbb{R}%
^{n+1}$ onto $T_{q}S^{n}=q^{\bot }$. In particular $\widetilde{\phi _{t}}$
fixes $\pm v\in S^{n}$. For $\tau =\widetilde{\phi _{t}}$ we will need
specifically the following standard facts:

$*$) If $u\in S^{n}$ is orthogonal to $v$, then $v\in T_{u}S^{n}$ and if $%
\tau \left( u\right) =\left( \cos \theta \right) u+\left( \sin \theta
\right) v$, then $\left( d\tau \right) _{u}v$ is a vector in $T_{\tau \left(
u\right) }S^{n}$ spanned by $u$ and $v$ of length $\cos \theta $.

$**$) There exists a positive constant $c$ such that $\left( d\tau \right)
_{\pm v}$ is a multiple $c^{\pm 1}$ of the identity map on $T_{\pm
v}S^{n}=v^{\bot }$.

Let $\Delta _{n}=\left\{ \left( p,p\right) \mid p\in S^{n}\right\} $ denote
the diagonal in $S^{n}\times S^{n}$. The map
\begin{eqnarray}
\psi :\mathcal{G}_{n+1}\rightarrow \left( S^{n}\times S^{n}\right)
\backslash \Delta _{n}\text{,\ \ \ \ \ }\psi \left( \left[ \gamma \right]
\right) =\left( \gamma \left( -\infty \right) ,\gamma \left( \infty \right)
\right)  \label{gss}
\end{eqnarray}
is a well-defined diffeomorphism. We denote by $\widehat{g}$ the induced
action of $g$ $\in G$ on $\left( S^{n}\times S^{n}\right) \backslash \Delta
_{n}$, that is $\widehat{g}\left( p,q\right) =\left( \widetilde{g}\left(
p\right) ,\widetilde{g}\left( q\right) \right) $. Given distinct points $%
p,q\in S^{n}$, let $T_{p,q}$ denote the reflection on $\Bbb{R}^{n+1}$ with
respect to the hyperplane orthogonal to $p-q$.

\begin{proposition}
If $\mathcal{G}_{n+1}$ is endowed with the metric $g_{1}$ and one considers
on $\left( S^{n}\times S^{n}\right) \backslash \Delta _{n}$ the
pseudo-Riemannian metric whose associated norm is
\begin{equation}
\left\| \left( x,y\right) \right\| _{\left( p,q\right) }=4\left\langle
T_{p,q}x,y\right\rangle /\left| q-p\right| ^{2}  \label{mss}
\end{equation}
for $x\in p^{\bot }$, $y\in q^{\bot }$, then the diffeomorphism $\psi $ of $%
\emph{(}$\emph{\ref{gss})} is an isometry.
\end{proposition}

\noindent \textbf{Proof.} Clearly $\psi $ is $G$-equivariant. Since the
metric $g_{1}$ on $\mathcal{G}_{n+1}$ is $G$-invariant, it is sufficient to
show that the metric (\ref{mss}) on $\left( S^{n}\times S^{n}\right)
\backslash \Delta _{n}$ is $G$-invariant as well and that $d\psi _{\left[
\gamma _{o}\right] }$ is a linear isometry.

Given distinct points $p_{\pm }\in S^{n}$, we show first that for any $g\in
G $ with $\widetilde{g}\left( e_{\pm 1}\right) =p_{\pm }$, $d\widehat{g}%
_{\left( -e_{1},e_{1}\right) }$ is a linear isometry. A straightforward
computation shows that the given metric on $\left( S^{n}\times S^{n}\right)
\backslash \Delta _{n}$ is invariant by the action of $SO_{n+1}$, since for
all $k$ in this group, $T_{k\left( p\right) ,k\left( q\right) }\circ
k=k\circ T_{p,q}$ for all $p,q\in S^{n}$, $p\ne q$. Hence we may suppose
without loss of generality that $p_{\pm }=\pm \left( \cos \theta \right)
e_{1}+\left( \sin \theta \right) e_{2}$ for some $\theta \in [0,\pi /2)$.
Now, any directly conformal transformation $\widetilde{g}$ as above may be
written as a composition $\tau ^{2}\circ \tau ^{1}\circ R$, where $R$ is a
rotation fixing $e_{1}$ and $\tau ^{1}$ and $\tau ^{2}$ are transvections
fixing $\left( -e_{1},e_{1}\right) $ and $\left( -e_{2},e_{2}\right) $,
respectively.

The assertion ($**$) above, with $v=e_{1}$ and $\tau =\tau ^{1}$, implies
that $d\widehat{\tau ^{1}}_{\left( -e_{1},e_{1}\right) }$ is a linear
isometry. Now we use the assertion ($*$) with $v=e_{2}$ and $u=e_{1}$ to see
that $d\widehat{\tau ^{2}}_{\left( -e_{1},e_{1}\right) }:e_{1}^{\bot }\times
e_{1}^{\bot }\rightarrow p_{-}^{\bot }\times p_{+}^{\bot }$ is a linear
isometry. Let $\lambda _{\pm }v+x_{\pm }\in T_{\pm u}S^{n}=u^{\bot }$, with $%
\lambda _{\pm }$ real numbers and $\left\langle x_{\pm },v\right\rangle =0$.
One computes
\begin{eqnarray}
\left\| \left( \lambda _{-}v+x_{-},\lambda _{+}v+x_{+}\right) \right\|
_{\left( -u,u\right) } &=&4\left( \lambda _{-}\lambda _{+}+\left\langle
x_{-},x_{+}\right\rangle \right) /\left| 2u\right| ^{2}  \label{tal} \\
&=&\left( \lambda _{-}\lambda _{+}+\left\langle x_{-},x_{+}\right\rangle
\right) \text{.}  \nonumber
\end{eqnarray}
On the other hand, call $d\tau _{\pm u}^{2}\left( v\right) =v_{\pm }$ and $%
d\tau _{\pm u}^{2}\left( x_{\pm }\right) =y_{\pm }$. Hence $\left| v_{\pm
}\right| =\cos \theta $. Since $d\tau _{\pm u}^{2}$ is conformal, $y_{\pm }$
is orthogonal to $v_{\pm }$ and has length $\left| x_{\pm }\right| \cos
\theta $. Also, $y_{\pm }$ is orthogonal to $u$, hence it is left fixed by $%
T_{p_{-},p_{+}}$. Therefore one computes
\[
\left\| \left( \lambda _{-}v_{-}+y_{-},\lambda _{+}v_{+}+y_{+}\right)
\right\| _{\left( -u,u\right) }=\frac{4\cos ^{2}\theta}{\left|
p_{-}-p_{+}\right| ^{2}}\left( \lambda _{-}\lambda _{+}+ \left\langle
x_{-},x_{+}\right\rangle \right) \text{,}
\]
which coincides with (\ref{tal}) since $\left| p_{-}-p_{+}\right| =2\cos
\theta $. This completes the proof that $d\widehat{g}_{\left(
-e_{1},e_{1}\right) }$ is a linear isometry. It remains only to show that $%
d\psi _{\left[ \gamma _{o}\right] }$ is a linear isometry.

We have that $\gamma _{o}\left( t\right) =\left( \cosh t,\sinh t,0\right)
\in \Bbb{R}^{n+2}$. Let $J$ be the Jacobi field along $\gamma _{o}$
orthogonal to $\gamma _{o}$ and satisfying $J\left( 0\right) =x$ and $%
J^{\prime }\left( 0\right) =y$, both in $T_{e_{0}}H$ orthogonal to $%
e_{1}=\gamma ^{\prime }\left( 0\right) $. We show next that
\[
d\psi _{\left[ \gamma _{o}\right] }L_{\gamma _{o}}\left( J\right) =\left(
x-y,x+y\right) \text{,}
\]
where $L_{\gamma _{o}}$ was defined in (\ref{el}). By invariance of $\psi $
by rotations it is sufficient to see that
\begin{eqnarray}
d\psi _{\left[ \gamma _{o}\right] }L_{\gamma _{o}}\left( J_{\pm }\right)
=\left( \pm e_{2},e_{2}\right) \text{,}  \label{psielpm}
\end{eqnarray}
where $J_{-}\left( 0\right) =0$, $J_{-}^{\prime }\left( 0\right) =e_{2}$, $%
J_{+}\left( 0\right) =e_{2}$ and $J_{+}^{\prime }\left( 0\right) =0$. Let
now
\[
A_{s}=\left(
\begin{array}{cc}
\cos s & -\sin s \\
\sin s & \cos s
\end{array}
\right) \text{ and }B_{s}=\left(
\begin{array}{ccc}
\cosh s & 0 & \sinh s \\
0 & 1 & 0 \\
\sinh s & 0 & \cosh s
\end{array}
\right) \text{.}
\]
The field $J_{-}$ is associated to the variation of $\gamma _{o}$
corresponding to the one parameter group of isometries $s\mapsto A_{s}^{-}=$
diag\thinspace $\left( 1,A_{s},I_{n-1}\right) $. One computes $%
A_{s}^{-}\left( \gamma _{o}\left( t\right) \right) =\left( \cosh t\right)
e_{0}+\sinh t\left( \left( \cos s\right) e_{1}+\left( \sin s\right)
e_{2}\right) \in H$. Hence
\begin{eqnarray*}
\left( A_{s}^{-}\circ \gamma _{o}\right) \left( \pm \infty \right)
&=&\lim_{t\rightarrow \pm \infty }\left( \tanh t\right) \left( \left( \cos
s\right) e_{1}+\left( \sin s\right) e_{2}\right) \\
&=&\pm \left( \cos s\right) e_{1}\pm \left( \sin s\right) e_{2}\text{,}
\end{eqnarray*}
whose derivative at $s=0$ is $\pm e_{2}$. Therefore $\left( \left.
d/ds\right| _{0}\right) \psi \left[ A_{s}^{-}\circ \gamma _{o}\right]
=\left( -e_{2},e_{2}\right) $. Using $B_{s}^{+}=$ diag\thinspace $\left(
B_{s},I_{n-1}\right) $ instead of $A_{s}^{-}$ one verifies the remaining
identity of (\ref{psielpm}). Finally, since $T_{-e_{1},e_{1}}$ clearly fixes
$x,y$, the norm (\ref{mss}) of $\left( x-y,x+y\right) $ at $\left(
-e_{1},e_{1}\right) $ is $4\left\langle x-y,x+y\right\rangle /\left|
2e_{1}\right| ^{2}=\left| x\right| ^{2}-\left| y\right| ^{2}$, which
coincides with the norm of $L_{\gamma _{o}}\left( J\right) $ by (\ref{nj}).
This shows that $d\psi _{\left[ \gamma _{o}\right] }$ is a linear isometry.
\hfill $\square $

\bigskip%

Let $\Bbb{O}$ denote the normed division algebra of the octonions and let $%
\Bbb{R}^{7}=$ Im\thinspace $\Bbb{O}$ endowed with its canonical cross
product $\times $. Let $j$ be the almost complex structure of $S^{6}$
defined by $j_{p}\left( x\right) =p\times x$ if $x\in T_{p}S^{6}=p^{\bot }$.
For $q\in S^{6}$, $q\ne p$, let $j_{p,q}$ be the linear operator on $%
T_{q}S^{6}=q^{\bot }$ defined by $j_{p,q}=T_{p,q}\circ j_{p}\circ T_{p,q}$.

\begin{proposition}
For all $x\in p^{\bot },y\in q^{\bot }$,
\[
J_{\left( p,q\right) }\left( x,y\right) =(j_{p}\left( x\right)
,j_{p,q}\left( y\right) )
\]
defines an orthogonal almost complex structure on $\left( S^{6}\times
S^{6}\right) \backslash \Delta _{n}$ with the metric above.
\end{proposition}

\noindent%
\textbf{Proof. }First we check that $J$ is an almost complex structure.
Indeed,
\[
\left\langle j_{p,q}\left( y\right) ,q\right\rangle =\left\langle
j_{p}T_{p,q}\left( y\right) ,T_{p,q}\left( q\right) \right\rangle
=\left\langle p\times T_{p,q}\left( y\right) ,p\right\rangle =0
\]
and $J^{2}=-\,$id holds as well, since $j_{p}^{2}=-$\thinspace id and $%
T_{p,q}^{2}=$ id. Finally, $J$ is orthogonal since both $j_{p}$ and $T_{p,q}$
are so. \hfill $\square $

\bigskip%

\noindent%
\textbf{Remarks. }a) By Proposition \ref{proptoni} there exists no proper
subgroup of $G$ acting transitively on $\mathcal{G}$ leaving $J$ invariant,
as it is the case of the analogous almost complex structure defined in \cite
{salvaimm} on the space of oriented lines of $\Bbb{R}^{7}$.

b) The structure $J$ is not integrable, since $\left( S^{6}\backslash
\left\{ p\right\} \right) \times \left\{ p\right\} $ is an almost complex
submanifold for any $p$, whose induced almost complex structure is $q\mapsto
j_{q}$, which is not integrable.

\bigskip

\noindent \textbf{Acknowledgment.} I would like to thank Eduardo Hulett for
his help and Antonio Di Scala for the statement and the idea of the proof of
Proposition \ref{proptoni}.



\noindent Marcos Salvai, \textsc{famaf\,-\,ciem}, Ciudad Universitaria, 5000
C\'ordoba, Argentina.

\smallskip

\noindent E-mail: salvai@mate.uncor.edu

\end{document}